\pdfoutput=1
\documentclass[12pt]{article}

\usepackage{amssymb,latexsym}

\usepackage{graphicx}
\usepackage{epstopdf}
\usepackage{enumerate}
\usepackage{pdfpages}
\usepackage{layout}
\usepackage[a4paper, total={6in, 8in}]{geometry}
\usepackage{polski}

\makeatother

\newtheorem{thm}{Theorem}[section]

\newtheorem{lem}[thm]{Lemma}

\newtheorem{prop}[thm]{Proposition}

\newcommand{\bd}{{\rm bd}}

\newcommand{\diam}{{\rm diam}}
\newcommand{\width}{{\rm width}}

\begin{document}

\baselineskip=18pt

\includepdf[pages=-]{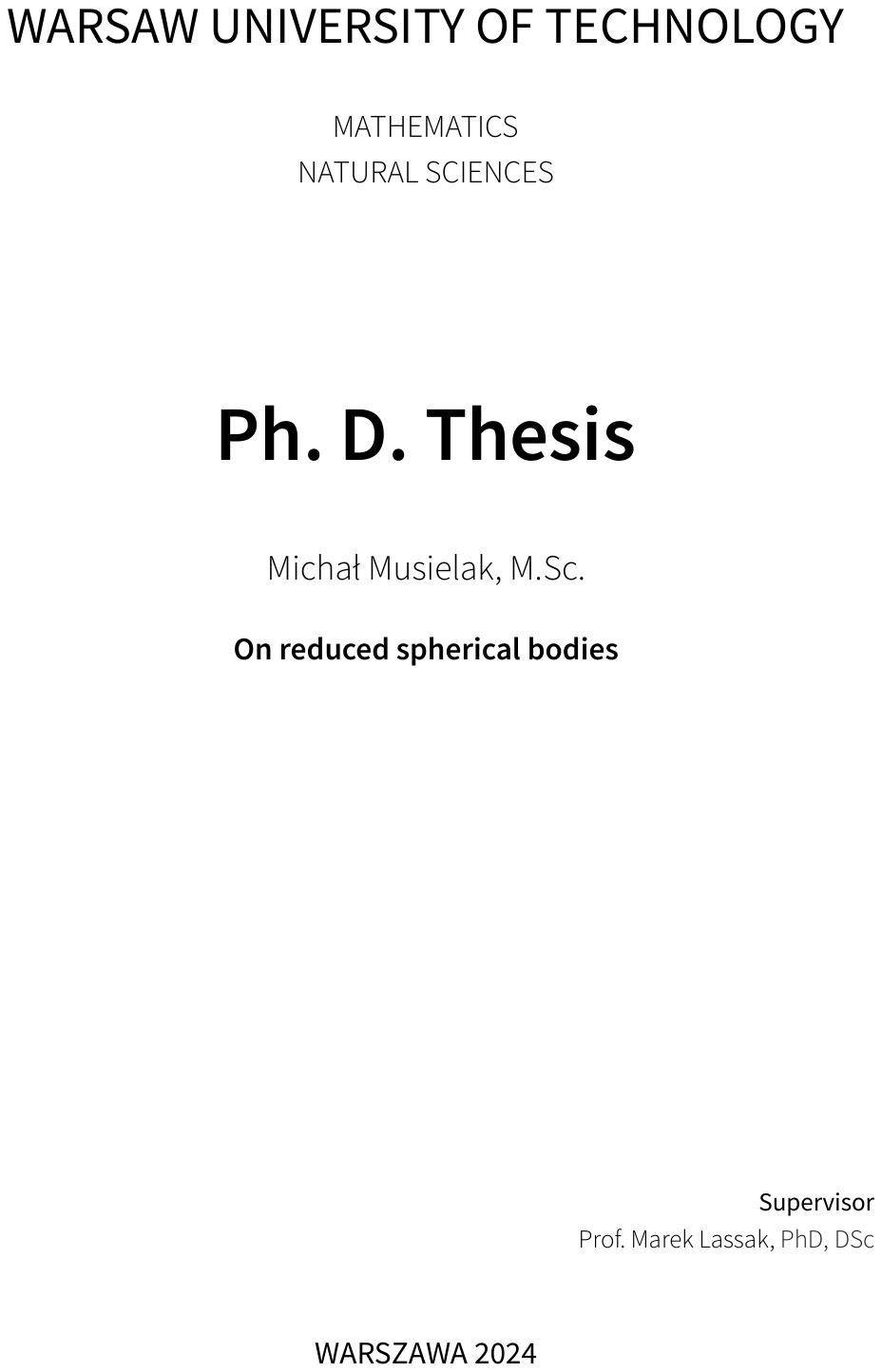}

\newpage
\thispagestyle{empty}
\

\newpage
\ 

\bigskip

\bigskip

\bigskip

\bigskip

\bigskip

\bigskip

\bigskip

\bigskip

\bigskip

{\Large \textbf{ACKNOWLEDGEMENTS}}

\bigskip

The author would like to thank Marek Lassak for his support and patience.

\newpage
\thispagestyle{empty}
\

\newpage

{\Large \textbf{ABSTRACT}}:

\medskip

This thesis consists of five papers about reduced spherical convex bodies and in particular spherical bodies of constant width on the $d$-dimensional sphere $S^d$.

In paper I we present some facts describing the shape of reduced bodies of thickness under $\frac{\pi}{2}$ on $S^2$. 
We also consider reduced bodies of thickness at least $\frac{\pi}{2}$, which appear to be of constant width.

Paper II focuses on bodies of constant width on $S^d$. 
We present the properties of these bodies and in particular we discuss conections between notions of constant width and  of constant diameter.

In paper III we estimate the diameter of a reduced convex body.

The main theme of paper IV is estimating the radius of the smallest disk that covers a reduced convex body on $S^2$.

The result of paper V is showing that every spherical reduced polygon $V$ is contained in a disk of radius equal to the thickness of this body centered at a boundary point of $V$.

\bigskip

\bigskip

{\Large \textbf{STRESZCZENIE}}:

\medskip

Niniejsza rozprawa sk\l ada si\k{e} z pi\k{e}ciu publikacji dotycz\k{a}cych sferycznych cia\l\ wypuk\l ych zredukowanych i w szczeg\'{o}lno\'{s}ci sferycznych cia\l\  o sta\l  ej szeroko\'{s}ci.

W publikacji I przedstawiamy pewne fakty dotycz\k{a}ce kszta\l tu cia\l\  zredukowanych o grubo\'{s}ci mniejszej ni\.z $\frac{\pi}{2}$ na $S^2$. 
Rozwa\.zamy r\'{o}wnie\.z cia\l a zredukowane o grubo\'{s}ci co najmniej $\frac{\pi}{2}$, kt\'{o}re okazuj\k{a} si\k{e} mie\'{c} sta\l \k{a} szeroko\'{s}\'{c}.

Publikacja  II skupia si\k{e} na cia\l ach o sta\l ej szeroko\'{s}ci na $S^d$. 
Przedstawiamy w\l asno\'{s}ci takich cia\l , a w szczeg\'{o}lno\'{s}ci omawiamy zwi\k{a}zek mi\k{e}dzy poj\k{e}ciami sta\l ej szeroko\'{s}ci i sta\l ej \'{s}rednicy.

W publikacji III szacujemy \'{s}rednic\k{e} zredukowanego cia\l a wypuk\l ego.

G\l \'{o}wnym tematem publikacji IV jest oszacowanie promienia najmniejszego dysku pokrywaj\k{a}cego zredukowane cia\l o wypuk\l e na $S^2$.

Rezultatem publikacji V jest wykazanie, \.ze ka\.zdy sferyczny zredukowany wielok\k{a}t $V$ jest zawarty w dysku o promieniu $\Delta (V)$ o \'{s}rodku w punkcie brzegowym $V$.

\bigskip

\bigskip

{\large \textbf{KEYWORDS}}: spherical convex body, spherical geometry, hemisphere, lune, width, constant width, thickness

\bigskip

{\large \textbf{SŁOWA KLUCZOWE}}: wypukłe ciało sferyczne, geometria sferyczna, hemisfera, luna, szerokość, stała szerokość, grubość

\newpage


\section*{Preface}

This thesis consists of a brief summary and the following papers.

\begin{description}
\item[I]  M. Lassak, M. Musielak, \textit{Reduced spherical convex bodies}, Bull. Pol. Ac. Math., \textbf{66} (2018), No 1, 87-97.
\item[II]  M. Lassak, M. Musielak, \textit{Spherical convex bodies of constant width}, Aequationes Math., \textbf{92} (2018), 627–640.
\item[III] M. Lassak, M. Musielak, \textit{Diameter of reduced spherical convex bodies}, Fasciculi Math., \textbf{61} (2018), 83–88.
\item[IV]  M. Musielak, \textit{Covering a reduced spherical body by a disk},  Ukr. Math. J., \textbf{72} (2021), 1613–1624.
\item[V] M. Musielak, \textit{Covering a reduced spherical polygon by a disk}, Rend. Circ. Mat. Palermo, II. Ser (2023). https://doi.org/10.1007/s12215-023-00899-z
\end{description}
\newpage

\tableofcontents

\newpage

\section{Introduction}

\subsection{Background}

The notion of a reduced convex bodies first appeared in 1987 in \cite{He} and was considered in $E^d$. 
Heil, the author of this paper, put forward the question of whether every reduced strictly convex body is of constant width. 
In the following years the subject was further developed by many authors: Dekster, Lassak, Martini, Groemer, Averkov and others. Some of important papers on that matter are: \cite{ChGr}, \cite{L1}, \cite{De2}, \cite{Sch}, \cite{MS}, \cite{Gr}, \cite{Av}, \cite{Gr2}.

\smallskip

In 2011 Lassak and Martini wrote the paper \cite{LM}, reviewing what was known about the reduced convex bodies and what open problems remained. 
An example of such a problem is the question of the existence of the reduced polytopes in $E^d$ for $d\ge 3$. 
It is solved in some particular cases, but in general the answer is still unknown.

\smallskip

There are some natural questions: how can we define reduced convex bodies in other geometries, which results can be transferred from the Euclidan spaces and if there are any new results substantially different from those in $E^d$.
In particular we can ask about the situation on the sphere $S^d$.
The first attemp to tackle this subject was made by Lassak in 2015 in \cite{L2}.
On the sphere we find the main problem caused by the lack of a notion of parallelity.
Therefore, we not only can not use methods of proving similar to $E^d$, but we have to actually find a new definitions leading us to the notion of reduced convex bodies.
Lassak has overcome that problem and presented such a definition, which allowed to begin research on reduced spherical convex bodies. 

\smallskip

One of the most important results from \cite{L2} is Theorem 1 describing the way of finding the width of convex bodies on the sphere. 
With that tool Lassak could come up with further interesting results.
Among those is a pointing some connections between the width of spherical convex body and its diameter.
Theorem 4 \cite{L2} says that \emph{through every extreme point $e$ of a reduced body $R\subset S^d$ a lune $L\supset R$ of thickness $\Delta(R)$ passes with $e$ as the center of one of the two $(d-1)$-dimensional hemispheres bounding $L$.} It leads to an open problem: is the same true for every boundary point of $R$?
We should also mention the fact proven in \cite{L2} that every smooth reduced convex body on the sphere is of constant width.

\smallskip

The paper \cite{L2} is the base for further research about reduced spherical convex bodies.

\subsection{Basic notions}

Consider the unit sphere $S^d$ in the $(d+1)$-dimensional Euclidean space $E^{d+1}$ for $d\geq 2$. 
By a {\it great circle} of $S^d$ we mean the intersection of $S^d$ with any two-dimensional subspace of $E^{d+1}$.
The common part of $S^d$ with any hyper-subspace of $E^{d+1}$ is called a {\it $(d-1)$-dimensional great sphere}.
Observe that each $1$-dimensional great sphere of $S^d$ is nothing else but a great circle.
By a pair of {\it antipodes} 
 of $S^d$ we understand a pair of points of $S^d$ such that one is a reflection of another with respect to the origin of $E^{d+1}$.

\smallskip

Clearly, if two different points $a, b \in S^d$ are not antipodes, there is exactly one great circle containing them.
As the {\it arc} $ab$ connecting those points $a$ and $b$ we define the shorter part of the great circle containing these points. 
The length of this arc is called the {\it spherical distance $|ab|$ of $a$ and $b$}. 
Moreover, we agree that the distance of coinciding points is $0$, and that of any pair of antipodes is $\pi$.

A {\it spherical ball $B_\rho(x)$ of radius $\rho \in (0, {\frac{\pi}{2}}]$}, or {\it a ball} in short is the set of points of $S^d$ at distances at most $\rho$ from a fixed point $x$, which is called the {\it center} of this ball. 
An {\it open ball} (a {\it sphere}) of radius $\rho$ is the set of all points of $S^d$ having distance smaller than (respectively, exactly) $\rho$ from a fixed point.
A spherical ball of radius $\frac{\pi}{2}$ is called a {\it hemisphere}. 
So it is the common part of $S^d$ and a closed half-space of $E^{d+1}$.
We denote by $H(m)$ the hemisphere with center $m$.
Two hemispheres with centers at a pair of antipodes are called {\it opposite}.

A {\it spherical $(d-1)$-dimensional ball of radius $\rho \in (0, {\frac{\pi}{2}}]$} is the set of points of a $(d-1)$-dimensional great sphere of $S^d$ which are at distances at most $\rho$ from a fixed point.
We call it the {\it center} of this ball.  
The $(d-1)$-dimensional balls of radius $\frac{\pi}{2}$ are called {\it $(d-1)$-dimensional hemispheres}, and {\it semicircles} for $d=2$.

A set $C \subset S^d$ is said to be {\it convex} if no pair of antipodes belongs to $C$ and if for every $a, b \in C$ we have $ab \subset C$.  
A closed convex set on $S^d$ with non-empty interior is called a {\it convex body}.

Since the intersection of every family of convex sets is also convex, for every set $A \subset S^d$ contained in an open hemisphere of $S^d$ there is the smallest convex set ${\rm conv} (A)$ containing $Q$. 
We call it {\it the convex hull of} $A$.

Let $C \subset S^d$ be a convex body. 
Let $Q \subset S^d$ be a convex body or a hemisphere.  
We say that $C$ {\it touches $Q$ from inside} if $C \subset Q$ and ${\rm bd} (C) \cap {\rm bd} (Q) \not = \emptyset$. 
We say that $C$ {\it touches $Q$ from outside} if $C \cap Q \not = \emptyset$ and ${\rm int} (C) \cap {\rm int} (Q) = \emptyset$. 
In both cases, points of ${\rm bd} (C) \cap {\rm bd} (Q)$ are called {\it points of touching}.
In the first case, if $Q$ is a hemisphere, we also say that $Q$ {\it supports} $C$, or {\it supports $C$ at $t$}, provided $t$ is a point of touching.
If at every boundary point of $C$ exactly one hemisphere supports $C$, we say that $C$ is {\it smooth}.
An extreme point of a convex body $C \subset S^2$ is a point for which the set $C \setminus \{e\}$ is convex. 
We denote the set of extreme points of $C$ by $E(C)$.

If hemispheres $G$ and $H$ of $S^d$ are different and not opposite, then $L = G \cap H$ is called {\it a lune} of $S^d$. 
This notion is considered in many books and papers (for instance, see \cite{VB}). 
The $(d-1)$-dimensional hemispheres bounding $L$ and contained in $G$ and $H$, respectively, are denoted by $G/H$  and $H/G$.

Observe that $(G/H) \cup (H/G)$ is the boundary of the lune $G \cap H$.  
Denote by $c_{G/H}$ and $c_{H/G}$ the centers of $G/H$ and $H/G$, respectively.
By {\it corners} of the lune $G \cap H$ we mean points of the set $(G/H) \cap (H/G)$. 
In particular, every lune on $S^2$ has two corners. 
They are antipodes. 

We define the {\it thickness $\Delta (L)$ of a lune} $L = G \cap H$ on $S^d$ as the spherical distance of the centers of the $(d-1)$-dimensional hemispheres $G/H$ and $H/G$ bounding $L$.

Compactness arguments show that for any hemisphere $K$ supporting a convex body $C \subset S^d$ there is at least one hemisphere $K^*$ supporting $C$ such that the lune $K \cap K^*$ is of the minimum thickness.
In other words, there is a ``narrowest" lune of the form $K \cap K'$ over all hemispheres $K'$ supporting $C$. 
The thickness of the lune $K \cap K^*$ is called {\it the width of $C$ determined by $K$.} 
We denote it by ${\rm width}_K (C)$. 

We define the {\it thickness} $\Delta (C)$ of a spherical convex body $C$ as the smallest width of $C$. 
This definition is analogous to the classical definition of thickness (called also minimal width) of a convex body in Euclidean space.

By {\it the relative interior} of a convex set $C \subset S^d$ we mean the interior of $C$ with respect to the smallest sphere $S^k \subset S^d$ that contains $C$. 

\smallskip

We define the polar set of a set $F\subset S^d$ as $F^o= \bigcap_{p\in F} H(p)$.
This notion is used in paper I and V.
In paper IV we define also $F^o_\rho$ as $\bigcap_{p\in F} B_{\rho}(p)$.
Clearly $F_{\pi / 2}^o= F^o$.

\smallskip

If $C \subset S^2$ is a spherical convex body and $X=H(x)$, $Y=H(x)$, $Z=Z(x)$ are different supporting hemispheres of $C$, then we write $\prec \!XYZ$ and say that $X, Y, Z$ {\it support $C$ in this order}, if $x,y,z$ are in this order on $\bd (F^\circ)$ when viewed from inside $C^{\circ}$.

\smallskip

In paper III we use a notion of a body of constant diameter. 
A body $W$ is of constant diameter $w$ if its diameter is equal $w$ and for every boundary point $p\in W$ there exists $p'\in W$ such that $|pp'|=w$

\subsection{Summary}

Paper \cite{L2} introduced most of the above notions.
This was a first and very important step in the study of spherical convex bodies.

\smallskip

This dissertation consists of papers that continue this research in different directions.
All of them are dedicated to the subject of reduced convex bodies on the sphere $S^d$.
They all represent the candidate's contribution to the development of this branch of science.

\smallskip

However it was not easy to make the progress at the beginning, since for a long time there was a problem with the proof of a fact analogous to Theorem 3 of \cite{L1}.

\smallskip

Finally, however, this problem was solved in paper I.
Theorem 3.1 from this work describes a situation analogous to that of the theorem mentioned in the last paragraph, but on a sphere.
This result made it possible to prove more facts analogous to those on the Euclidean plane.
Altogether, it can be estimated that the candidate's contribution to this work is equal to $50\%$, since all the results were obtained in teamwork with the co-author.

\smallskip

Paper III is, in a sense, a complement to paper I, as it also presents analogous facts to those in the paper \cite{L1}, but on a sphere.
As in paper I, during the research work the candidate and the co-author collaborated to achieve all the results.
Their contributions were more or less equal, so we can estimate the candidate's contribution at $50\%$.

\smallskip

The candidate was the sole author of papers IV and V, and therefore contributed fully to them.
These papers were further work on obtaining some results analogous to those of the plane. 
Paper IV is devoted to such a result similar to the plane one from paper \cite{L3}, and paper V presents such a result similar to that from paper \cite{Fa}.

\smallskip

Compared to the above papers, paper II deals with a slightly more specific topic, as it deals with the special case of reduced spherical bodies, i.e. spherical bodies of constant width, but in this case in any dimension.
Since we considered here the situation on the sphere $S^d$ for all $d \ge 2$, in the candidate's opinion this was both the most difficult and the most rewarding research paper of all those mentioned in this dissertation. 
Once again, the candidate and co-author collaborated and contributed equally to this work, which means that the candidate contributed $50\%$.

\medskip

As can be seen, the subject matter of all the work is very consistent. 
The candidate was concerned with a very specific area of convex geometry.
The aim of this dissertation is to present some results on reduced bodies on the sphere when this topic has already been studied in depth on the Euclidean plane and when two papers \cite{L2} and \cite{L4} have already given the basic properties of reduced bodies on the sphere.
The research work was also devoted to transferring some results on reduced convex bodies from the Euclidean plane to the sphere.

The following sections provide a more detailed discussion of these results.

\newpage


\section{Results in paper I}

The main aim of this paper is to find similar results to those established by Lassak for the Euclidean plane in \cite{L1} and compare the situation in the plane and on the sphere.

\smallskip

Since there is no parallelity on the sphere, we must use a completely different approach than the one in the plane, but it appears that spherical reduced convex bodies of thickness less that $\frac{\pi}{2}$ have similar properties as those on the plane.

\smallskip

We present two following facts describing the shape of the boundary of spherical reduced convex bodies.

\begin{thm} \label{main}
Let $R \subset S^2$ be a reduced spherical body with $\Delta (R) < \frac{\pi}{2}$.
Let $M_1$ and $M_2$ be supporting hemispheres of $R$ such that ${\rm width}_{M_1} (R) = \Delta (R) = {\rm width}_{M_2} (R)$ and ${\rm width}_{M} (R) > \Delta (R)$ for every hemisphere $M$ satisfying $\prec \!M_1MM_2$.
Consider lunes $L_1 = M_1 \cap M_1^*$ and $L_2 = M_2 \cap M_2^*$.  
Then arcs $a_1a_2$ and $b_1b_2$ are in $\bd (R)$, where $a_i$ is the center of $M_i/M_i^*$ and $b_i$ is the center of $M_i^*/M_i$ for $i=1,2$.
Moreover, $|a_1a_2| = |b_1b_2|$.
\end{thm}

Proof of this theorem consists of two parts.
In the first one we show that arcs $a_1a_2$ and $b_1b_2$ lay in the boundary of $R$.
In both cases we imagine the opposite, i.e. that there is an extreme point of $R$ between $a_1$ and $a_2$, and between $b_1$ and $b_2$. 
Then we reach a contradiction.
In the second part considering some two congruent triangles leads us to the equality $|a_1a_2| = |b_1b_2|$. 
\begin{flushright}
$\centerdot$
\end{flushright}

\smallskip

The next theorem presents the fact about a hemisphere supporting a reduced body.

\begin{thm} \label{segment}
Let $R$ be a reduced body with $\Delta (R) < \frac {\pi}{2}$.
Assume that $M$ is a supporting hemisphere of $R$ such that the intersection of $\bd(M)$ with $\bd (R)$ is a non-degenerated arc  $x_1x_2$. 
Then ${\rm width}_M (R) = \Delta (R)$, and the center of $M/M^*$ belongs to $x_1x_2$.
\end{thm}

The main idea of the proof of this theorem is assuming that the lune $M\cap M^*$ has the width greater than $\Delta (R)$ and next obtaining a contradiction.

\medskip

The following proposition leads us to the next theorem, which asserts that there is no reduced spherical convex body of width of at least $\frac{\pi}{2}$, that is not of constant width

\begin{prop}\label{smooth1}
Each reduced spherical body of thickness over $\frac{\pi}{2}$ is smooth.
\end{prop}

\begin{thm}\label{constant} 
Every reduced spherical convex body $R$ such that $\Delta (R) \geq \frac{\pi}{2}$ is a spherical body of constant width $\Delta (R)$. 
\end{thm}

For $\Delta(R)> \frac{\pi}{2}$ the thesis is an obvious consequence of the previous proposition.
For $\Delta(R) =\frac{\pi}{2}$ let us imagine the opposite, i.e. that there exist such hemispheres $K$ and $K^*$ that the width of the lune $K\cap K^*$ is greater than $\Delta (R)$.
Then we show the existence of an extreme point $e$ of $R$, for which the lune $K\cap H(e)$ is of width $\Delta (R)$.
But it is a contradiction, because from we know that $\Delta(K\cap K^*)=\Delta (K\cap H(e))$.

\bigskip

The fifth part of the paper is dedicated to a spherical convex body of constant width. 
We prove the following theorem.

\begin{thm}\label{strictly1} 
A reduced spherical convex body of thickness below $\frac{\pi}{2}$ is of constant width if and only if it is strictly convex.
\end{thm}

If a reduced spherical convex body of thickness below $\frac{\pi}{2}$ is not of constant width, then an evident consequence of Theorem \ref{main} is that it is not strictly convex.
The opposite implication requies a much more advanced proof.
To accomplish this we assume that there exists a not strictly convex body $W$ of constant width.
Then we point out the hemisphere $H$ and a point $c$, such that $c\notin H\cap H^*$ and this contradiction ends the proof.

\medskip

The last important result of the paper is a partial answer to a question asked by Lassak in \cite{L2}.

\begin{thm} \label{question}
Let $W$ be a spherical body of constant width.
For every boundary point $p$ of $W$ there exists a lune $L \supset W$ fulfilling $\Delta (L) = \Delta(W)$ such that $p$ is the center of one of the semicircles bounding $L$.
\end{thm}

The only non-trivial part of the proof is a case when $p$ is not an extreme point.
Therefore due to Theorem \ref{strictly1} the thesis holds true for bodies of constant width less than $\frac{\pi}{2}$.
For bodies of constant width equal $\frac{\pi}{2}$ we show that an expected lune is the intersection of the hemisphere centered at $p$ and the hemisphere supporting $W$ at $p$.
In the last case of bodies of constant width over $\frac{\pi}{2}$ we construct the ball of a radius $\Delta(W)- \frac{\pi}{2}$ whose center coincides with the center of the hemisphere supporting $W$ at $p$.
Then we show that a subset of the boundary of this ball is also a boundary of $W$, with ends the proof.

\newpage


\section{Results in paper II}

Paper II  focuses on spherical bodies of constant width, but considered on the sphere $S^d$ for any $d\ge 2$, not only $d=2$.

\smallskip

Among a few lemmas presented in the paper one should be mentioned here, due to its importance for further proofs.

\begin{lem}\label{convexhull}
Let $o \in S^d$ and $0 < \mu < \frac{\pi}{2}$.
For every $x \in S^d$ at distance $\frac{\pi}{2}$ from $o$ denote by $x'$ the point of the arc $ox$ at distance $\mu$ from $x$. 
Consider two points $x_1,x_2$ at distance $\frac{\pi}{2}$ from $o$ such that $|x_1x_2| < \pi - \mu$. Then for every $x \in x_1x_2$ we have 

$$B_\mu (x')  \subset {\rm conv} (B_\mu(x_1') \cup B_\mu(x_2')).$$
\end{lem} 

The second part of the paper opens with the theorem useful in the rest of the paper.

\begin{thm}\label{touching ball}  At every boundary point $p$ of a body $W \subset S^d$ of constant width $w > \pi/2$ we can inscribe a unique ball $B_{w- \pi/2}(p')$ touching $W$ from inside at $p$. 
What is more, $p'$ belongs to the arc connecting $p$ with the center of the unique hemisphere supporting $W$ at $p$, and $|pp'|=w-\frac{\pi}{2}$.
\end{thm}

In the proof we simply point out a ball mentioned in the theorem.
If $p$ is an extreme point of $W$, then we consider a lune $K\cap M$ such that $p$ is  the center of $K/M$.
We denote by $m$ the center of $M$.
Then we show that the ball centered in $m$ and with the radius $w- \frac{\pi}{2}$ is a ball mentioned in our theorem.
The case in which $p$ is not an extreme point is in some way similar, but a bit more complicated.
We know that in this case $p$ belongs to such a convex hull of some extreme points of $W$, which is a part of $\bd(W)$.
Finally, applying Lemma \ref{convexhull} and the previous case, we show that the thesis holds true for every point of this convex hull, so in particular for the point $p$.

\medskip

The next facts are a generalizations of the results obtained in \cite{LMu} for $d=2$. 
Here we present facts for any integer $d \geq 3$, however methods of proving are different than for the case $d=2$.  

\begin{prop}\label{smooth2} 
Every spherical body of constant width larger than $\frac{\pi}{2}$ (and more general, every reduced body of thickness larger than $\frac{\pi}{2}$) of $S^d$ is smooth. 
\end{prop}

\begin{thm} \label{strictly2} 
Every spherical convex body of constant width smaller than $\frac{\pi}{2}$ on $S^d$ is strictly convex. 
\end{thm}

In a short proof we imagine the opposite and using this assumption we find an arc in the boundary of our body which has a length greater that $\Delta(W)$.
This leads to the conclusion that $\diam (W)> \Delta (W)$, which contradicts the fact proven earlier in this paper that $\diam (W)\le \frac{\pi}{2}$.
This entails the equality $\diam (W) = \Delta (W)$.

\medskip

The following theorem gives a partial answer to the question asked by Lassak in \cite{L2}.
The theorem shows that the answer is positive in the case of bodies of constant width.

\begin{thm} \label{center} 
For every body $W \subset S^d$ of constant width $w$ and every $p \in {\rm bd} (W)$ there exists a lune $L \supset W$ fulfilling $\Delta (L) = w$ with $p$ as the center of one of the two $(d-1)$-dimensional hemispheres bounding this lune.
\end{thm}

Since the case when $p$ is an extreme point is an easy one, the thesis is obvious for $w< \frac{\pi}{2}$ since then every boundary point is extreme.
If $w=\frac{\pi}{2}$ and $p$ is not an extreme point, then we show that our lune is $H(p)\cap G$, where $G$ is a hemisphere supporting $W$ at $p$.
The most difficult case is when $w>\frac{\pi}{2}$ and $p$ is not an extreme point.
We construct a set denoted by $U$ and take the ball $B$ touching $W$ from inside.
Then in a very long and technical way we show that the intersection of $\bd (W)$ with $Q$ is equal to the intersection of $\bd(W)$ with $U$. 
Using that fact we point out the lune asserted in the theorem.

\medskip

We also show the following fact for bodies of constant width.

\begin{thm} \label{diam=w}
If $W \subset S^d$ is a body of constant width $w$, then ${\rm diam} (W)=w$.
\end{thm}

For $w \le \frac{\pi}{2}$ the thesis is a consequence of an analogous fact already proven in \cite{L2}.
For $w > \frac{\pi}{2}$ we consider an arc in the boundary of $W$ whose length is equal to $\diam (W)$. 
Then we take the lune $L$ for which one of the endpoints of this arc is the center of a $(d-1)$-dimensional hemisphere bounding $L$, and such that $\Delta(L)=\Delta(W)$.
We show that the width of this lune is not less than the length of our arc and therefore $w\ge \diam (W)$.
Since we also know that $w \le \diam (W)$, we arrive at our thesis.

\medskip

In the final part of paper II we consider the bodies of a constant diameter.
The last theorem in this paper presents the relationship between these bodies and bodies of constant width.

\begin{thm} \label{iff} 
Every spherical convex body $W \subset S^d$ of constant width $w$ is of constant diameter $w$.
Every spherical convex body $W \subset S^d$ of constant diameter $w \ge \frac{\pi}{2}$ is of constant width $w$. 
\end{thm}

The first statement of this theorem is an easy consequence of a fact already proven in this paper. 
The second statement is equivalent to the statement that for every hemisphere $K$ supporting $W$ we have $\width_K(W)=w$.
We consider separately the cases $w> \frac{\pi}{2}$ and $w= \frac{\pi}{2}$.
In both of them we assume the opposite and in a technical way obtain a contradicition.

\medskip

The paper ends with the question about the generalization of the second part of the last theorem: is it true that also every spherical body of constant diameter $w < \frac{\pi}{2}$ is a body of constant width $w$?
To the best of the author’s knowledge  it is still an open problem.

\newpage


\section{Results in paper III}

In \cite{L1} Lassak proved the inequality $\frac{\diam (R)}{\Delta(R)}\le \sqrt{2}$, where $R$ is any reduced convex body in the Euclidan plane.
There is a natural idea to show an analogous fact on the sphere $S^2$.
This is the theme of paper III.

For estimating the diameter of a spherical convex reduced body, we apply the following simple lemma:

\begin{lem} \label{diamE}
For every convex body $C \subset S^2$ of diameter at most $\frac{\pi}{2}$ we have\\ ${\rm diam} (E(C)) = {\rm diam}(C)$.
\end{lem}

Thanks to this lemma it is sufficient to focus only on extreme points of a reduced body.

\smallskip

The main theorem of the paper presents an estimate of the diameter:

\begin{thm} \label{ineq} 
For every reduced spherical body $R \subset S^2$ with $\Delta (R) < \frac{\pi}{2}$ we have ${\rm diam}(R) \leq {\rm arc cos} (\cos^2 \Delta (R))$.  
This value is attained if and only if $R$ is the quarter of disk of radius $\Delta(R)$.
If $\Delta (R) \geq \frac{\pi}{2}$, then ${\rm diam}(R) = \Delta (R)$.
\end{thm}

For proving the first part of the theorem, due to Lemma \ref{diamE}, it is sufficient to show that any two extreme points $e_1,e_2$ of $R$ are in a distance at least ${\rm arc cos} (\cos^2 \Delta (R))$.
Both those points are the centers of some semicircles bounding two lunes $L_1$ and $L_2$ containing $R$.
Let $f$ be such a point of an arc connecting the centers of two semicircles bounding $L_1$ such that the triangle $e_1fe_2$ is right.
We estimate the length of the legs of this triangle and therefore we are able to estimate its hypotenuse, which leads to the demanded value ${\rm arc cos} (\cos^2 \Delta (R))$.

\smallskip

The consequence of the last theorem is the following proposition, stating the connection between the size of a reduced body and its diameter.

\begin{prop} \label{precise}
Let $R \subset S^2$ be a reduced body.
Then ${\rm diam} (R) < \frac{\pi}{2}$ if and only if $\Delta (R) < \frac{\pi}{2}$. 
Moreover, ${\rm diam} (R) = \frac{\pi}{2}$ if and only if $\Delta (R) = \frac{\pi}{2}$. 
\end{prop}

\newpage


\section{Results in paper IV}

The paper gives an estimate for the radius of the smallest disk covering a spherical reduced convex body.
The same problem for Euclidean plane was solved by Lassak in \cite{L3}.
However, we are not able to apply the same method of proving, due to the lack of the notion of parallelity.
Therefore we apply a new, different approach.

The first step is the same as in the plane: we use the spherical Helly theorem, which in the spherical case says that if any three points of a convex body can be covered by a disk of a fixed radius, then the whole body can be covered by a disk of this radius.
As a consequence, it is sufficient to focus on three fixed points of our body.

\smallskip

The main idea is to start with solving the problem for a few simple cases and then reduce the general situation to these cases.
Three lemmas from the first part of the paper shows the value of the radius of the covering disk of particular figures: a quarter of a disk, a Realeaux triangle and a equilateral triangle.

\smallskip

Then we consider the case of a body of constant width. 
We apply the result of Dekster from \cite{De1} in the form useful for us: every spherical convex body of constant width at most $\frac{2\pi}{3}$ can be covered by a disk of radius $\arcsin \left( \frac{2\sqrt{3}}{3} \cdot \sin \frac{\Delta(W)}{2}\right)$.

An example of the Realeaux triangle 
 shows that this estimate can not be improved if $W$ is of constant width at most $\frac{\pi}{2}$.
But the following theorem presents the improvement in the case when $\Delta(W)$ is at least $\frac{\pi}{2}$.

\begin{thm}\label{const}
Every spherical body $W$ of constant width at least $\frac{\pi}{2}$ is contained in a disk of radius $\Delta(W) + \arcsin \left( \frac{2\sqrt{3}}{3} \cdot \cos \frac{\Delta(W)}{2}\right) - \frac{\pi}{2}$.\end{thm}

In the proof we use the notion $F^\oplus$, which is defined as $\left\{ p: F\subset H(p)\right\}$ for any set $F$ on the sphere.
It is easy the show that $F^\oplus$ is nothing else than the polar set $F^o$.  
However, it is more convenient to use the approach with $F^\oplus$. 

It is easy to show that if $W$ is a body of constant width, then $W^\oplus$ is the body of constant width $\frac{\pi}{2}-\Delta (W)$.
We apply this observation in the proof of our theorem.
It is convienent to look at the body $W^\oplus$ for which the radius of the covering disk can be estimated thanks to Dekster's result.
We can then calculate the growth of the radius for the body $W$.

\smallskip

The last part of the paper presents our main estimate. After proving two technical lemmas, we move into the main theorem of the paper:

\begin{thm}\label{second}
Every reduced spherical body $R$ of thickness at most $\frac{\pi}{2}$ is contained in a disk of radius $\rho = \arc\tan \left( \sqrt{2} \cdot \rm{tan} \frac{\Delta(R)}{2}\right)$.
\end{thm}

As mentioned earlier, it is sufficient to show that every three points od $R$ are contained in a disk of radius $\rho = \arc\tan \left( \sqrt{2} \cdot \tan \frac{\Delta(R)}{2}\right)$.
In order to show that we vary the position of our points in such a way, that the minimal radius of a disk covering them never decreases.
In some cases after moving the points we obtain 
a quarter of a disk, which ends the proof.
Otherwise, in the main case we obtain the situation, when our three points are edges of an isosceles triangle.
For this triangle we can find the connections between the values of its angles, sides and heights.
The obtained formula for the tangent of the radius of the disk circumscribed about the triangle is complicated, but using one of the previously proven lemmas we check that the greatest value of $\tan \rho$ is $\arc\tan \left( \sqrt{2} \cdot \rm{tan} \frac{\Delta(R)}{2}\right)$, which ends the proof.

\smallskip

Our estimate can not be improved, since the value from the theorem is attained for the quarter of the disk.

\newpage


\section{Results in paper V}

Fabi\'nska in \cite{Fa} proved that every reduced polygon in the Euclidan plane can be covered by the ball centered at a boundary point of the polygon whose radius is equal to the thickness of the polygon.
In paper V we show the analogon of this result for a reduced polygon on the sphere.
We apply here some of the results of the paper \cite{L4} which is dedicated to reduced spherical polygons.

\smallskip

We begin the paper with a few useful lemmas leading us to the corollary, which describes the shape of a polar set of a reduced polygon on the sphere.
Using this corollary we are able to prove the main theorem of this paper.

\begin{thm}
Every spherical reduced polygon $V$ is contained in a disk of radius $\Delta (V)$ centered at a boundary point of $V$.
\end{thm}

\smallskip

The main concept of the proof can be described in the following way.
We observe that the thesis of this theorem is equivalent to the statement that $V^\circ_{\Delta(V)}$ has a common point with the boundary of $V$.
We imagine the opposite and then using the previous corollary we obtain a contradiction. 
To obtain this contradiction we consider any arc $uv \subset V$ of length $\diam (V)$ and then we show that we can find a point $x\in V$ for which $|ux|>|uv|$. 
The knowledge about the shape of a polar set of a reduced polygon is substantial here.

\newpage

\

\bigskip

\bigskip

\bigskip

\bigskip

\bigskip

\bigskip

\bigskip

\bigskip

\bigskip

\bigskip

\bigskip

\bigskip

\bigskip

\bigskip

\bigskip

\bigskip

\bigskip

\bigskip

\bigskip

\bigskip

\bigskip

\bigskip

\bigskip

\bigskip

\begin{center}
{\Huge \textbf{Attachments}}
\end{center}

\newpage
\begin{center}
\

\bigskip

\

\bigskip

\

\bigskip

\bigskip

\bigskip

\bigskip

{\huge PAPER I}

\bigskip

\bigskip

\bigskip

\bigskip

{\huge  \textit{Reduced spherical convex bodies}}

\bigskip

\bigskip
 
{\Large M. Lassak, M. Musielak}

\bigskip

\bigskip

{\large Bull. Pol. Ac. Math. \textbf{66} (2018), No 1, 87-97.}
\end{center}

\newpage

\includepdf[pages=-,scale=.96,pagecommand={}]{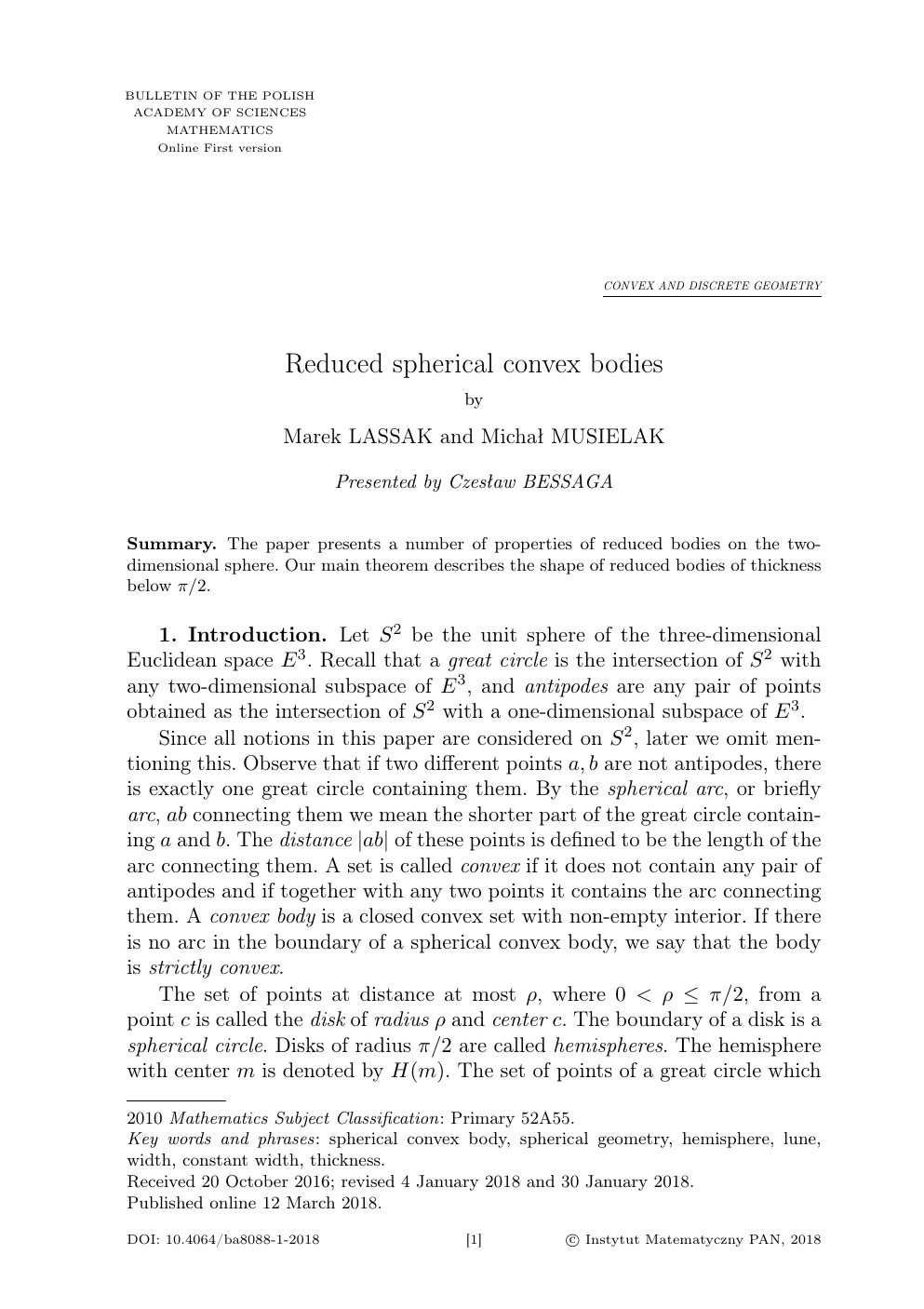}

\newpage

\begin{center}

\

\bigskip

\

\bigskip

\

\bigskip

\bigskip

\bigskip

\bigskip

{\huge PAPER II}

\bigskip

\bigskip

\bigskip

\bigskip

{\huge  \textit{Spherical convex bodies of constant width}}

\bigskip

\bigskip
 
{\Large M. Lassak, M. Musielak}

\bigskip

\bigskip

{\large Aequationes Math. \textbf{92} (2018), 627-640.}
\end{center}

\newpage

\includepdf[pages=-,scale=.90,pagecommand={}]{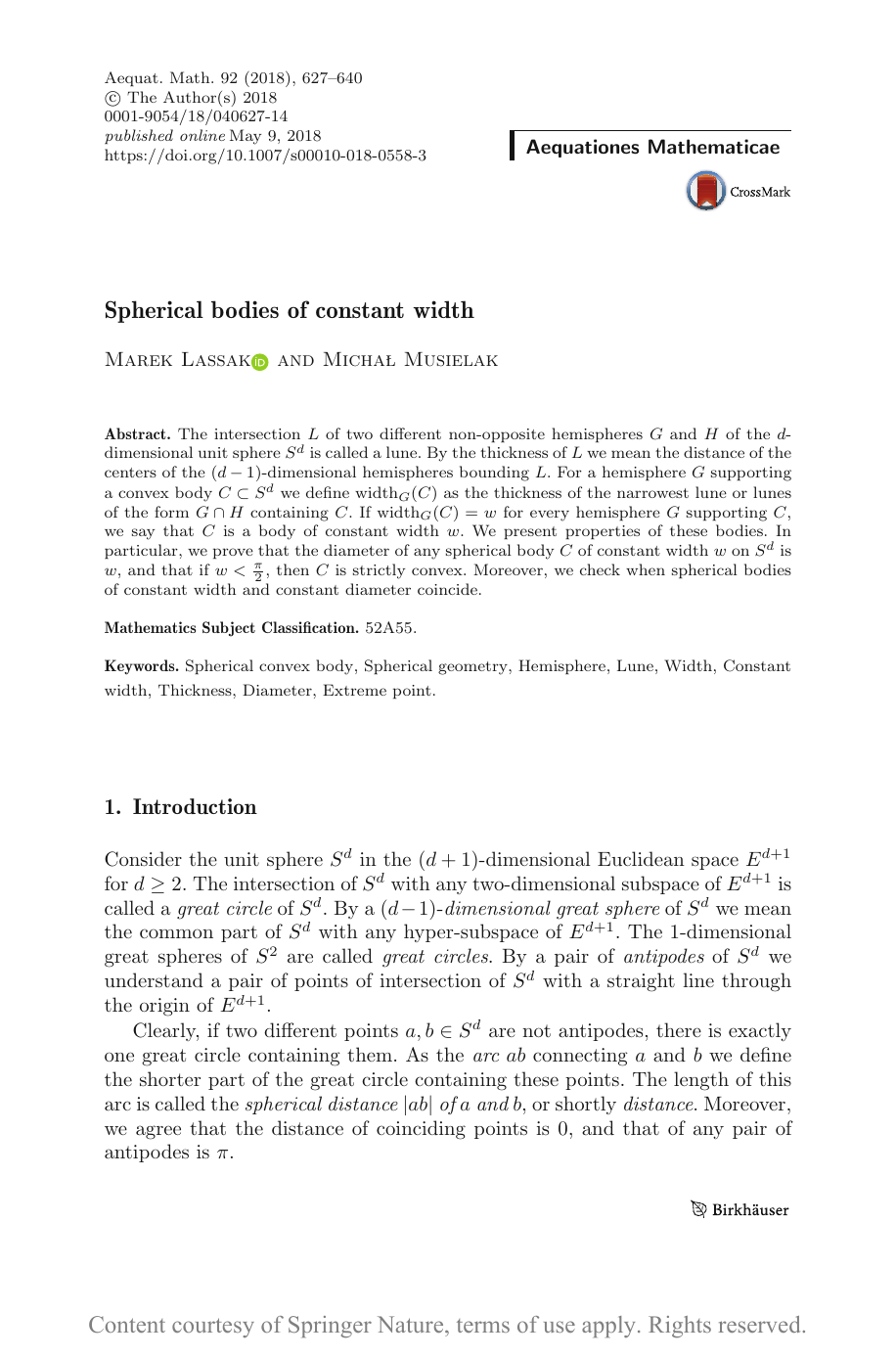}

\newpage

\begin{center}

\

\bigskip

\

\bigskip

\

\bigskip

\bigskip

\bigskip

\bigskip

{\huge PAPER III}

\bigskip

\bigskip

\bigskip

\bigskip

{\huge  \textit{Diameter of reduced spherical convex bodies}}

\bigskip

\bigskip
 
{\Large M. Lassak, M. Musielak}

\bigskip

\bigskip

{\large Fasciculi Math. \textbf{61} (2018), 83-88.}
\end{center}

\bigskip

\bigskip

\bigskip

\textit{The copyright for this publication is held by the journal 
Fasciculi Mathematici. Citation of this publication has been made possible by the courtesy of the editors of this journal, who have given their written permission.}

\newpage

\includepdf[pages=-,scale=.96,pagecommand={}]{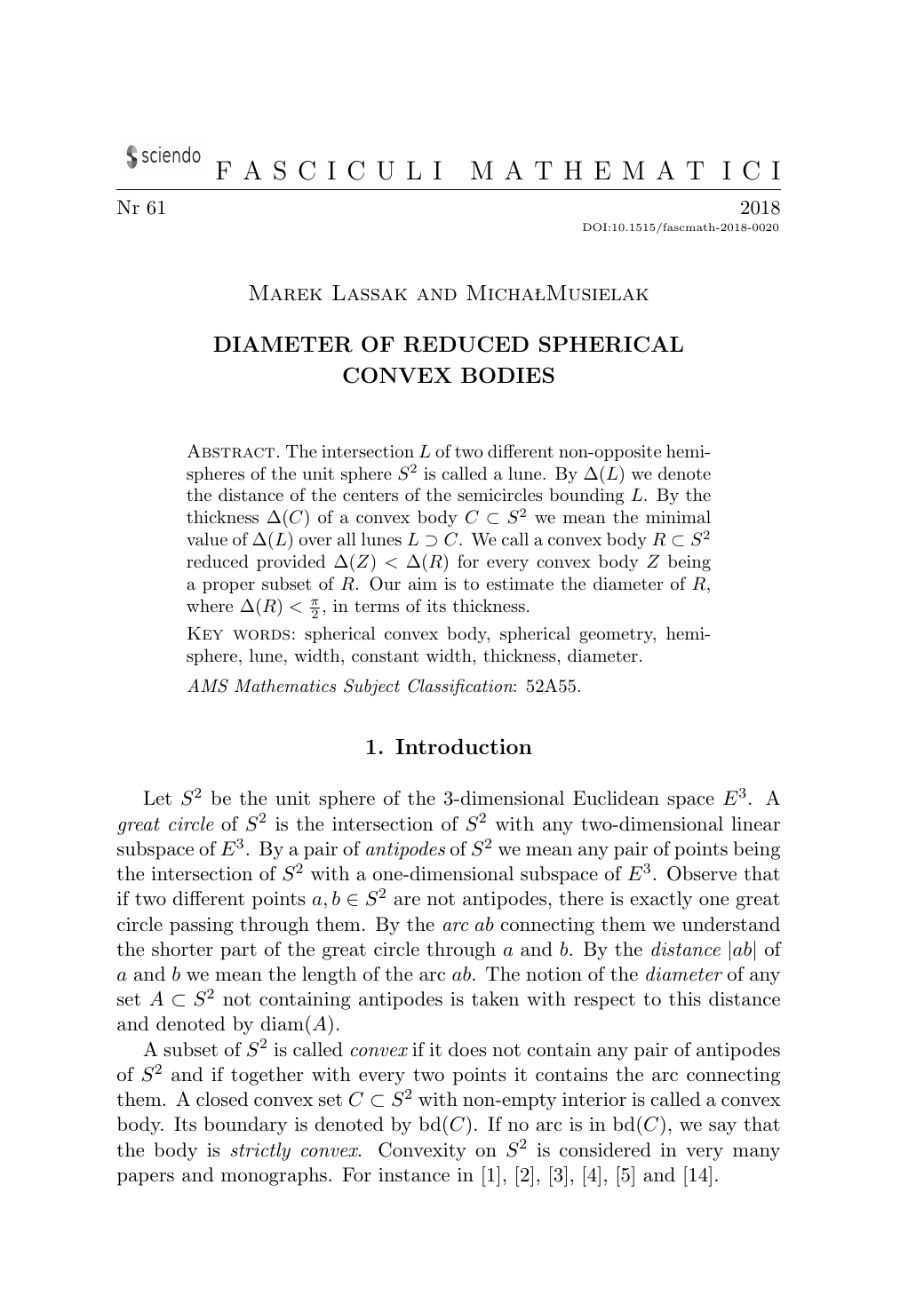}

\newpage

\begin{center}

\

\bigskip

\

\bigskip

\

\bigskip

\bigskip

\bigskip

\bigskip

{\huge PAPER IV}

\bigskip

\bigskip

\bigskip

\bigskip

{\huge  \textit{Covering a reduced spherical body by a disk}}

\bigskip

\bigskip
 
{\Large M. Musielak}

\bigskip

\bigskip

{\large Ukr. Math. J. \textbf{72}, 1613-1624 (2021).}
\end{center}

\newpage

\includepdf[pages=-,scale=.96,pagecommand={}]{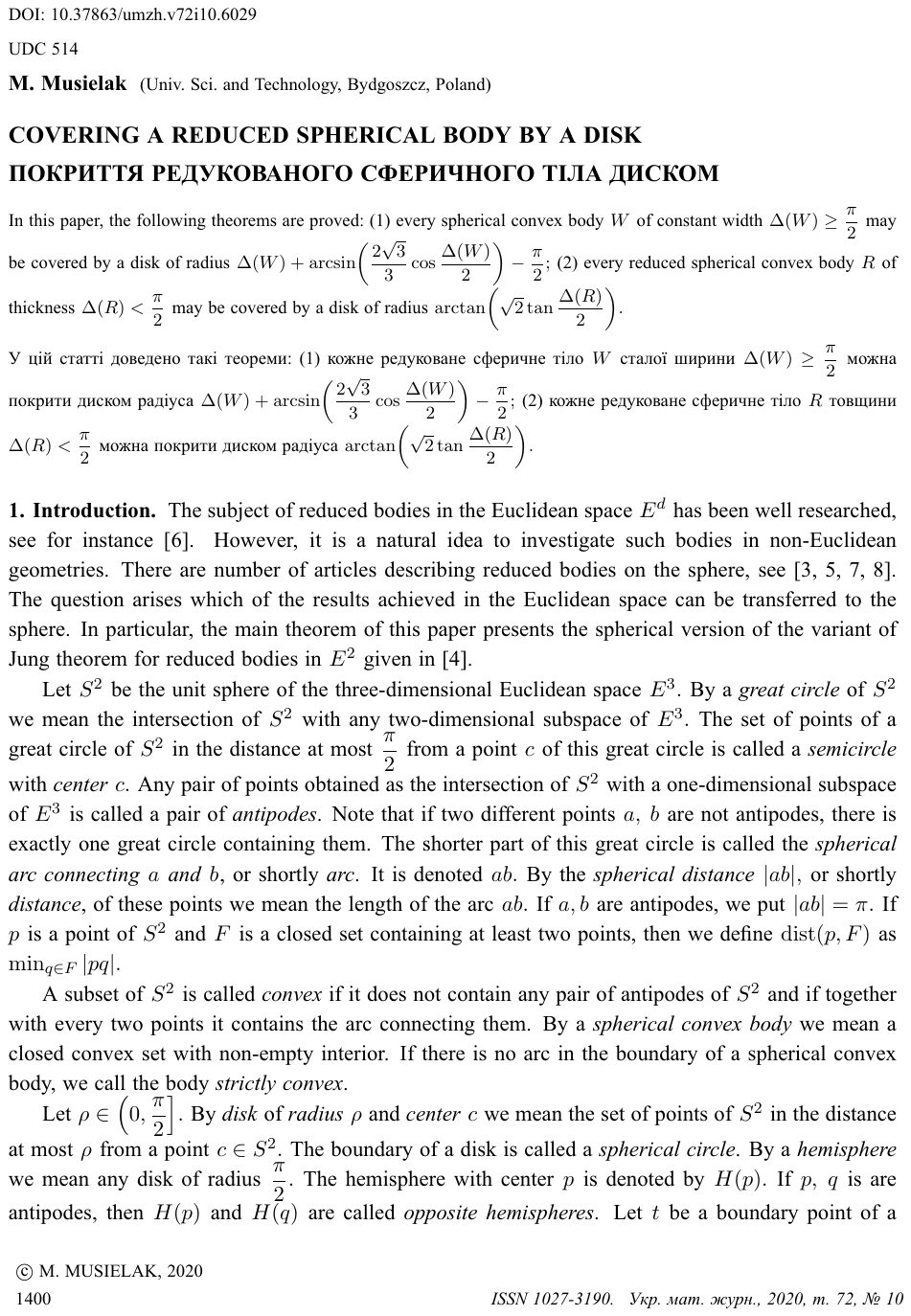}

\newpage

\begin{center}

\

\bigskip

\

\bigskip

\

\bigskip

\bigskip

\bigskip

\bigskip

{\huge PAPER V}

\bigskip

\bigskip

\bigskip

\bigskip

{\huge  \textit{Covering a reduced spherical polygon by a disk}}

\bigskip

\bigskip
 
{\Large M. Musielak}

\bigskip

\bigskip

{\large Rend. Circ. Mat. Palermo, II. Ser (2023).}
\end{center}

\newpage

\includepdf[pages=-,scale=.94,pagecommand={}]{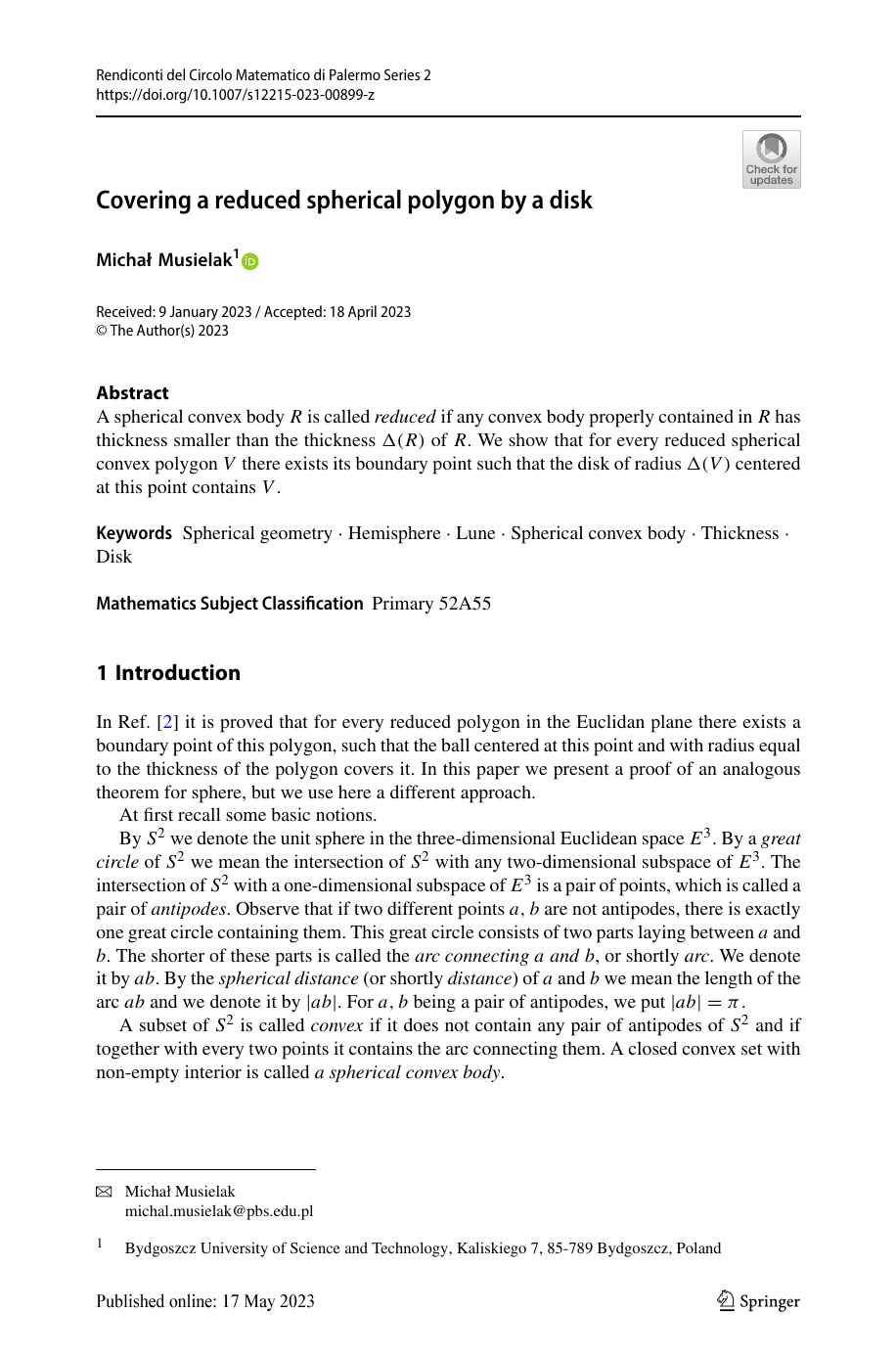}


\newpage

\begin{thebibliography}{15}

\bibitem{Av} G. Averkov, \textit{On planar convex bodies of given Minkowskian thickness and least possible area}, Arch. Math. \textbf{84} (2005)
183--192.


\bibitem{ChGr} G. D. Chakerian, H. Groemer, {\it Convex bodies of constant width}, Convexity and its applications, 49--96,
Birkhauser, Basel, 1983.
  

\bibitem {De1} B. V. Dekster, {\it The Jung theorem for spherical and hyperbolic spaces}, Acta Math. Hungar.,  {\bf 67 (4)} (1995), 315--331

\bibitem {De2} B.V. Dekster, \textit{On reduced convex bodies}, Israel J. Math. \textbf{56} (1986) 247--256.

\bibitem{Fa}  E. Fabi\'{n}ska, \textit{Covering a reduced polygon by a disk}, Demonstratio Math., \textbf{42} (2009), No. 4, 851--856


\bibitem{Gr2} H.Groemer, \textit{Extremal convex sets} Monatsh. Math., \textbf{96}, 29--39 (1983).

\bibitem{Gr} H. Groemer, \textit{On complete convex bodies}, Geom. Dedicata, \textbf{20} (1986) 319-–334.


\bibitem{He} E. Heil, \textit{Kleinste konvexe K\"{o}rper gegebener Dicke}, \#453, Fachbereich Mathematik, Technische Hochschule Darmstad, 1987.





\bibitem{L1}   
M. Lassak, \textit{Reduced convex bodies in the plane},  Israel J. Math., {\bf 70} (1990), 365--379. 
 
 \bibitem{L3} M. Lassak, {\it On the smallest disk containing a planar reduced convex body}, Arch. Math., \textbf{80} (2003), 553--560.
 
\bibitem{L2} M. Lassak, {\it Width of spherical convex bodies}, Aequationes Math., 
 \textbf{89} (2015), no.  3, 555--567.  
 
 \bibitem{L4} M. Lassak, \textit{Reduced spherical polygons},  Colloq. Math., {\bf 138} (2015), 205--216. 

\bibitem{LM} M. Lassak, H. Martini, {\it Reduced convex bodies in Euclidean space -- a survey}, Expositiones Math., {\bf 29} (2011), 204--219. 

\bibitem{LMu} M. Lassak, M. Musielak, {\it Reduced spherical convex bodies}, Bull. Pol. Ac. Math., \textbf{66} (2018), 87--97.

\bibitem{LMu2} M. Lassak, M. Musielak, {\it Spherical bodies of constant width}, Aequationes Math., \textbf{92} (2018), 627--640.

\bibitem{LMu3} M. Lassak, M. Musielak, \textit{Diameter of reduced spherical convex bodies}, Fasciculi Math., \textbf{61} (2018), 83–88.

\bibitem{Le} K. Leichtweiss, {\it Curves of constant width in the non-Euclidean geometry}, Abh. Math. Sem. Univ. Hamburg, {\bf 75} (2005), 257--284.

\bibitem{MS}  H. Martini, K.J. Swanepoel, \textit{The geometry of Minkowski spaces—a survey}, Part II, Expositiones Math., \textbf{22} (2004) 93--144.

  

\bibitem{Mu}  M. Musielak, \textit{Covering a reduced spherical body by a disk},  Ukr. Math. J., \textbf{72} (2021), 1613--1624.

\bibitem{Mu2} M. Musielak, \textit{Covering a reduced spherical polygon by a disk}, Rend. Circ. Mat. Palermo, II. Ser (2023)



\bibitem{Sch} R. Schneider, \textit{Convex bodies: the Brunn–Minkowski Theory}, Cambridge University Press, Cambridge, 1993.


\bibitem {VB} 
G. Van Brummelen, {\it Heavenly mathematics. The forgotten art of spherical trigonometry.}, Princeton University Press (Princeton, 2013). 


\end{thebibliography}
\end{document}